\documentclass [a4paper,12pt] {amsart}
\usepackage [latin1]{inputenc}
\usepackage[all]{xy}
\usepackage{latexsym}  
\usepackage{amsmath}
\usepackage{amsfonts}
\usepackage{amssymb}
\usepackage{graphics}
\usepackage{graphicx}
\newtheorem{teorema}{Theorem}[section]
\newtheorem{lemma}[teorema]{Lemma}
\newtheorem{propos}[teorema]{Proposition}
\newtheorem{corol}[teorema]{Corollary}
\newtheorem{ex}{Example}[section]
\newtheorem{rem}{Remark}[section]
\newtheorem{defin}[teorema]{Definition}
\def\bt{\begin{teorema}}
\def\et{\end{teorema}}
\def\bp{\begin{propos}}
\def\ep{\end{propos}}
\def\bl{\begin{lemma}}
\def\el{\end{lemma}}
\def\bc{\begin{corol}}
\def\ec{\end{corol}}
\def\br{\begin{rem}\rm}
\def\er{\end{rem}}
\def\bex{\begin{ex}\rm}
\def\eex{\end{ex}}
\def\bd{\begin{defin}}
\def\ed{\end{defin}}
\def\demo{\par\noindent{\bf Proof.\ }}
\def\enddemo{\ $\Box$\par\vskip.6truecm}
\def\R{{\mathbb R}}   \def\a {\alpha} \def\b {\beta}
\def\N{{\mathbb N}}     \def\d {\delta} 
\def\C{{\mathbb C}}      \def\l{\lambda}
                 \def\p{\partial}
                  \def\s{\sigma}

 \def\oli{\overline}

\def\O{\Omega}

\begin {document}
\title[Cohomology of semi $1$-coronae]{Cohomology of semi $1$-coronae and extension of analytic subsets}
\author[Alberto Saracco, Giuseppe Tomassini]{Alberto Saracco, \\ Giuseppe Tomassini}
               \address{Scuola Normale Superiore, Piazza dei Cavalieri, 7 - I-56126 Pisa, Italy}
              \email{a.saracco@sns.it, g.tomassini@sns.it}
             \date{\today}
             \keywords{q-pseudoconvexity $\cdot$ cohomology $\cdot$ bump lemma $\cdot$ extension}
              \subjclass[2000]{Primary 32D15 Secondary 32F10, 32L10}
\maketitle
\tableofcontents
\section{Introduction and notations}
Let $X$ be a (connected and reduced) complex space. We
recall that $X$ is said to be {\it strongly} $q$-{\it pseudoconvex} in the sense of
Andreotti-Grauert~\cite{AG} if there exists a compact subset $K$ and a smooth function
$\varphi:X\to\R$, $\varphi\ge 0$, which is strongly $q$-plurisubharmonic on $X\setminus
K$ and such that 
\begin{itemize}
\item[a)] for every $b>\max\limits_K\,\varphi$ the subset
$$ 
B_b=\{x\in X:\varphi(x)<b\}  
$$  
is relatively compact in $X$.  
\end{itemize}
If $K=\emptyset$, $X$ is said to be $q$-{\it complete}. We remark that, for a space, being $1$-complete is equivalent to being Stein.
Replacing the condition a) by 
\begin{itemize}
\item[a')] for every $\min\limits_X,\varphi<a<\min\limits_K\,\varphi$
and every $b>\max\limits_K\,\varphi$ the subset
$$ 
C_{a,b}=\{x\in X:a<\varphi(x)<b\}  
$$
is relatively compact in $X$,
\end{itemize}
we obtain the notion of $q$-{\it corona} (see~\cite{AG}, \cite{AT}).

A $q$-corona is said to be {\it complete} whenever $K=\emptyset$.

The cohomology of coherent sheaves defined on $q$-coronae was studied in \cite{AT}. 

In this paper we will focus on {\it semi $q$-coronae}, which are defined as follows. 

Consider a strongly $q$-pseudoconvex space (or, more generally, a $q$-corona) $X$, and a smooth function $\varphi:X\to\R$ displaying the
$q$-pseudoconvexity of $X$. Let $C_{a,b}$ be a $q$-corona of $X$ and let $h:X\to\R$ be a
pluriharmonic function (i.e.\ locally the real part of a holomorphic function) such that: $K\cap\{h=0\}=\emptyset$, $\{h=0\}$ and the boundaries of $B_b$ and $B_a$ are $\C$-transversal. A connected component $C^+_{a,b}$ of $C_{a,b}\setminus\{h=0\}$ is, by definition, a {\it semi} $q$-{\it corona}.

At the origin of the interest for domains whose boundary contains a Levi flat part there is an extension theorem for CR-functions proved in \cite{LT} (see also \cite{La}, \cite{LaP}, \cite{St}). 

A more general type of semi $q$-corona is obtained by replacing the zero set of $h$ with with a Levi flat hypersurface $H$.

In both cases the semi $q$-coronae are differences $B^+_b\setminus{\overline B^+}_a$
where $B^+_b$, $B^+_a$ are strongly $q$-pseudoconvex spaces. Indeed, the function
$\psi=-\log h^2$ (respectively $\psi=-\log \delta_H(z)$, where $\delta_H(z)$ is the
distance of $z$ from $H$) is plurisubharmonic in $W\setminus\{h=0\}$ (respectively in $W\setminus H$) where $W$ is a  neighbourhood of $B_b\cap \{h=0\}$ (respectively of $B_b\cap H$). Let $\chi:\R\to\R$ be an
increasing convex function such that $\chi\circ\varphi>\psi$ on a neighbourhood of
$B_b\setminus W$. The function $\Phi=\sup\,(\chi\circ\varphi,\psi)+\varphi$ is an exhaustion function 
for $B_b\setminus \{h=0\}$ (respectively for $B_b\setminus H$) and it is strongly $q$-plurisubharmonic in $B_b\setminus(\{h=0\}\cup K)$ (respectively in 
$B_b\setminus H\cup K)$.

Results on the cohomology of coherent sheaves on semi $q$-coronae were obtained by the authors in \cite{ST} under the hypothesis that they are defined on the larger set $B_b$. 

The aim of this paper is to give a generalization for coherent sheaves $\mathcal F$ defined only on the semi $q$-corona. For the sake of simplicity we restrict ourselves to the case of smooth semi $1$-coronae.

Following Andreotti-Grauert (see \cite{AG}), given a semi $1$-corona 
$$
C^+_{a,b}=C_{a,b}\cap\{h>0\},$$
where $h$ is pluriharmonic, and a coherent sheaf $\mathcal F$ on $C^+_{a,b}$ we consider the strongly plurisubharmonic functions $P_\varepsilon(z)=\varepsilon\vert z\vert^2-h$, $\varepsilon>0$, and an exhaustion of $C^+_{a,b}$ by the following relatively compact domains 
$$
C^+_\varepsilon=\left\{z\in\C^n:P_\varepsilon(z)< 0\right\}\bigcap \oli C_{a+\varepsilon,b-\varepsilon}.
$$
 The idea is to  prove for the domains $ C^+_\varepsilon$  a bump lemma and an approximation theorem as in the classical case of coronae. Here the situation is more complicated because of the presence of a non-empty pseudoconvex-pseudoconcave part in the boundary of each $C^+_\varepsilon$, In order to circumvent this difficulty, we work with the closed sets $\oli C^+_\varepsilon$ using in a crucial way a regularity result on the $\oli\p$-equation due to Laurent-Thi\`ebaut and Leiterer  (see Section \ref{CORISO}). This  enables us to prove the following results: assume that ${\rm depth}\,\mathcal F_z\ge 3$ for $z$ near to the pseudoconcave part of the boundary of $C^+_{a,b}$; then\vspace{.3cm}
\begin{enumerate}
\item[1)] if $\varepsilon$ is sufficiently small and $\varepsilon'<\varepsilon$ is near $\varepsilon$ 
$$
H^1(\oli C^+_{\varepsilon'},\mathcal F)\simeq H^1(\oli C^+_{\varepsilon},\mathcal F) 
$$
\item[2)] the cohomology spaces $H^1(\oli C^+_{\varepsilon},\mathcal F)$ are finite dimensional.
\end{enumerate}
\vspace{.3cm}
(see Lemma \ref{BUMP1}, Lemma \ref{ONTO5} and Proposition \ref{ONTO51}). 

Thus the function
$$
d(\varepsilon)={\rm dim}_\C\,H^1(\oli C^+_{\varepsilon},\mathcal F)
$$
is piecewise constant, but, in general, it could have frequently a \lq\lq jump-discontinuity\rq\rq\ and it could happen that $d(\varepsilon)\to+\infty$ (see Remark \ref{JUMP}).
Neverthless, the isomorphism $1)$ allows us to prove in the last section:
\begin{enumerate}
	\item the fact that Oka-Cartan-Serre Theorem $A$ holds in semi $1$-coronae for sheaves which satisfy the condition of Theorem \ref{ISO} (see Theorem \ref{OKS});
	\item an extension theorem for analytic subsets (see Corollary \ref{OKA2}).
\end{enumerate}
 
It is worth noticing that an extension theorem for codimension one analytic subsets of a non-singular semi $1$-corona was proved in \cite{ST} and for higher codimensions, using different methods based on Harvey-Lawson's theorem \cite{HL}, by Della Sala and the first author in \cite{DS}. 

\section{Remarks on the proofs of theorems in \cite{ST}}
Let $X$ be a complex space. For every coherent sheaf $\mathcal F$ on $X$ and every subset $A$ of $X$ we set 
$$
p(A;\mathcal F)=\inf\limits_{x\in A}{\rm depth}(\mathcal F_x)
$$ 
$$
p(A)=p(A;\mathcal O).
$$
Let $C=C_{a,b}$ be a $q$-corona of $X$. All the results in \cite{ST} on finite and/or vanishing cohomology results for $q$-coronae and semi $q$-coronae are obtained using Andreotti-Grauert methods. They consist of two main points 
\begin{itemize}
\item[i)]the bump lemma;
\item[ii)]for every corona  $C_{a',b'}\Subset C$ there exist a corona $C_{a'+\varepsilon,b'+\varepsilon}\Subset C $, $\varepsilon>0$ such that the homomorphism
$$
H^r(C_{a'-\varepsilon,b'+\varepsilon},\mathcal F)\longrightarrow H^r(C_{a',b'},\mathcal F)
$$
is bijective for $q\leq r\leq p(C;\mathcal F).$
\end{itemize}
As a matter of fact the method of proof shows that the condition on the depth is needed only in $C_{a,a'}$ i.e.\ the homomorphism  
$$
H^r(C_{a'-\varepsilon,b'+\varepsilon},\mathcal F)\longrightarrow H^r(C',\mathcal F)
$$
is bijective for $q\leq r\leq p(C_{a,a'};\mathcal F). $ 

Let $X$ be a strongly $q$-pseudoconvex space (respectively $X\subset\C^n$ be a strongly $q$-pseudoconvex open set) and $H=\left\{h=0\right\}$ where $h$ is pluriharmonic in $X$ (respectively $H$ Levi-flat), and $C=C_{a,b}=B_b\setminus{\overline B}_a$ a $q$-corona. We can suppose that $B_b \setminus H$ has two connected components, $B^+$ and $B^-$, and we define $C^+=B^+\cap C$, $C^-=B^-\cap C$.

From the above remark we derive the following improvements of Theorem 1, Corollary 2 and Theorem 3 in \cite{ST}. 
\bt\label{Ac}
Let $\mathcal F\in {\sf Coh}(B_b)$. Then the image of the homomorphism 
$$
H^r(\oli{B}^+,\mathcal F)\oplus H^r(\oli{C},\mathcal F)\longrightarrow H^r(\oli{C}^+,\mathcal F)
$$
(all closures are taken in $B_b$), defined by $(\xi\oplus\eta)\mapsto\xi_{|\oli{C}^+}-\eta_{|\oli{C}^+}$ has finite codimension provided that $q-1\le
r\le p(\oli B_a;\mathcal F)-q-2$.
\et
\bc\label{cB}
If $K\cap H=\emptyset$, under the same assumption of Theorem~\ref{Ac}
$$
\dim_\C\,H^r(\oli C^+,\mathcal F)<\infty
$$
for $q\le r\le p({\overline B_a};\mathcal F)-q-2$.
\ec
\bt\label{cC}
If $\oli B_+$ is a $q$-complete space, then
$$
H^r(\oli C,\mathcal F)\stackrel{\sim}{\rightarrow} H^r(\oli C^+,\mathcal F)
$$
for $q\le r\le p(\overline B_a;\mathcal F)-q-2$, and the homomorphism
\begin{equation}\label{eqA}
H^{q-1}(\oli B^+,\mathcal F)\oplus H^{q-1}(\oli C,\mathcal F)\longrightarrow H^{q-1}(\oli C^+,\mathcal F)
\end{equation}
is surjective for $p({\overline B_a};\mathcal F)\ge 2q+1$.

If $\oli B^+$ is a $1$-complete space and $p({\overline B_a};\mathcal F)\ge 3$, then
$$
H^0(\oli B^+,\mathcal F)\stackrel{\sim}{\rightarrow}H^0(\oli C^+,\mathcal F).
$$
\et
This implies the following. Let $C_1=B_{b_1}\smallsetminus {\oli B}_{a_1}\Subset C_2=B_{c_2}\smallsetminus {\oli B}_{a_2}$. Then
$$
H^r(\overline C^+_1,\mathcal M_{\{x\}}\mathcal F)\stackrel{\sim}{\rightarrow}H^r(\overline C^+_2,\mathcal F),
$$
for $q\leq r\leq p({\overline B_{a_1}};\mathcal F)$.

In particular, if $x\in C_2\smallsetminus \overline B_{a_1}$ and $\mathcal M_{\{x\}}$ denotes the sheaf of ideals of $\{x\}$, then
$$
H^r(\overline C^+_2,\mathcal M_{\{x\}}\mathcal F)\stackrel{\sim}{\rightarrow}H^r(\overline C^+_1,\mathcal F),
$$
for $q\leq r\leq p({\overline B_{a_1}};\mathcal F)$.

\section{An isomorphism theorem for semi $1$-coronae}\label{CORISO}
Our aim is to give a generalization of the above results for sheaves defined only on the semi $q$-coronae, i.e.\ for the case when the \lq\lq hole\rq\rq\ is real. For the sake of simplicity will consider only complete $1$-coronae in $\C^n$ with $n\ge 3$. So we consider connected $1$-coronae of the form
$$
C=\{z\in\C^n\ :\ 0< \varphi(z)<1\}\Subset\C^n,
$$
where $\varphi:\C^n\to\R$ is a smooth strongly plurisubharmonic function in a Stein neighborhood $U$ of $\{0\le\varphi\le1\}$, $d\varphi\neq 0$ on $\varphi=0,1$. Let $h$ be a pluriharmonic function on $U$ and $H$ the zero set of $h$. We assume that $H$ is smooth and transversal to the hypersurfaces $\left\{\varphi=0\right\}$, $\left\{\varphi=1\right\}$, that $U\smallsetminus H$ has two connected components $U^\pm$ and $h>0$ on $U^+$. For $0<a<b<1$ we set
$$
B_{b}=\{z\in U:\varphi< b\},\>B^+_{b}=B_{b}\cap U^+,
$$
$$
C_{a,b}=\left(B_{b}\smallsetminus \overline B_a\right),\> C^+_{a,b}=C_{a,b}\cap U^+.
$$
Let $P_\varepsilon(z)=\varepsilon\vert z\vert^2-h$; then there is $\varepsilon_0$ such that for $\varepsilon\in (0,\varepsilon_0)$ the hypersurfaces $\{\varphi=\varepsilon\}$, $\{\varphi=1-\varepsilon\}$ meet $\{P_\varepsilon=0\}$ transversally. Finally we define the following subsets (which are locally $1$-{\em convex}, $1$-{\em concave}, see \cite{LTL} and Remark \ref{cvxcnc} below) 
$$
\oli C^+_\varepsilon=\left\{z\in\C^n:P_\varepsilon(z)\le 0\right\}\bigcap \oli C_{\varepsilon,1-\varepsilon}.
$$
We want to prove the following 
\bt\label{ISO}
Let $C^+$ be a semi $1$-corona in $\C^n$.
Then for every $\varepsilon\in(0,\varepsilon_0)$ there exists $\oli\varepsilon\in[0,\varepsilon)$ such that for every $\mathcal F\in {\sf Coh}(C^+)$  satisfying
$$
 \oli {\left\{z\in C^+:{\rm depth}(\mathcal F_z)<3\right\}}\bigcap B_{\varepsilon_0}=\emptyset. 
$$
and every $\ \varepsilon'\in (\oli\varepsilon,\varepsilon)$ the homomorphism 
$$
H^1(C^+_{\varepsilon'},\mathcal F)\longrightarrow H^1(\oli C^+_\varepsilon,\mathcal F)
$$
is an isomorphism. 
\et
The main ingredients for the proof are the bump lemma and a density theorem as in Andreotti-Grauert \cite{AG}. In order to treat points belonging to the pseudoconvex-pseudoconcave part of the boundary we work with closed bumps using the following result due to Laurent-Thi\`ebaut and Leiterer (see \cite[Proposition 7.5]{LTL}):
\bp\label{LTL} Let $D\Subset\C^n$ be a $1$-concave, $1$-convex domain of order $1$ of special type, and suppose that $n\geq3$. If $f$ is a continuous ($n,r$)-form in some neighborhood $U_{\oli D}$ of $\oli D$, $1\leq r\leq n-2$, such that $\oli\partial f=0$ in $U_{\oli D}$, then there exists a form $u\in\bigcap_{\varepsilon>0} C^{1/2-\varepsilon}_{n,r-1}(\oli D)$ such that $\oli\partial u=f$ in $D$.
\ep
\br\label{cvxcnc} Proposition 7.5 in \cite{LTL} is much more general, but we state it this way, since the semi $1$-coronae we consider are locally $1$-concave, $1$-convex domain of order $1$ of special type, i.e.\ they are locally biholomorpic to the set-difference of two convex domains.\er
The proof of Theorem \ref{ISO} is a consequence of several intermediate results.
\subsection{Bump lemma: surjectivity of cohomology}
With the same notations as above let $\oli D=\oli C^+_\varepsilon$, $0<\varepsilon<\varepsilon_0$ where $\varepsilon_0<b$ is so chosen that for all $\varepsilon\in (0,\varepsilon_0)$ the hypersurfaces $\{\varphi=\varepsilon\}$, $\{\varphi=1-\varepsilon\}$ are $\C$-transversal to $\{P_\varepsilon=0\}$. Let $\Gamma_1$, $\Gamma_2$ be respectively the pseudoconvex and the pseudoconcave part of the boundary ${\rm b}D$ of $\oli D$. Thus ${\rm b}D=\oli \Gamma_1\cup\oli\Gamma_2$ and $\oli\Gamma_2$ is contained in the smooth hypersurface $\{\varphi=\varepsilon\}$.
\bl[bump lemma]\label{BUMP1}There exists a finite open covering $\mathcal{U}$ of ${\rm b}D$, $\mathcal U=\left\{\oli U_j\right\}_{1\le j\le
m}$, and compact subsets $\oli D_1,\ldots,D_m$ of $C^+$ such that
\begin{itemize}
\item[(i)] $\oli D=\oli D_0\subset\oli D_1\subset\cdots\subset \oli D_m $;\\
\item[(ii)] $\oli D\subset D_m$;\\
\item[(iii)] $\oli D_j\setminus \oli D_{j-1}\subset \oli U_j$ for $1\le j\le m$;\\
\item[(iv)] if $\mathcal F\in {\sf Coh}(C^+)$ then 
$$
H^r(\oli U_j\cap \oli D_k,\mathcal{F})=0
$$
for every $j,k$ and $1\le r\le p(\oli D;\mathcal F)-2$. 
\end{itemize} 
Moreover, the family of the coverings $\mathcal{U}$ as
above is cofinal in the family of all finite coverings of
${\rm b}D$.
\el
\demo
If $z^0\in\Gamma_1\cup\Gamma_2$ i.e.\ $z^0$ is a point of pseudoconvexity or pseudoconcavity we argue as in the proof of the classical Andreotti-Grauert bump lemma.

Assume that $z^0\in\oli\Gamma_1\cap\oli\Gamma_2$. There exists a a sufficiently small closed ball $\oli B$ of positive radius, centered at $z^0$ and a biholomorphism on $\Phi :\oli B\to \Phi(\oli B)$  which trasnsform $\oli B\cap \{\varphi\ge \varepsilon\}$ and $\oli B\cap\{P_\varepsilon\ge 0\}$ respectively in a strictly concave and strictly convex set. We may also assume that $\mathcal F_{|\oli B}$ has a homological resolution
\begin{equation}\label{Hom}
0\ \to\ \mathcal O^{p_k}\ \to\ \cdots\ \to\ \mathcal O^{p_0}\ \to\ \mathcal F\ \to\ 0
\end{equation}
with $n-k\ge 3$. Choose  a smooth function $\varrho\in C^{\infty}_0(B)$ such that $\varrho\ge 0$ and $\varrho (z^0)\neq 0$ and
a positive number $\lambda$ such that the closed domains
$$
\oli B_1=\left\{\varphi-\varepsilon-\lambda\rho \le 0\right\}\cap\oli B,\>\oli B_2=\left\{P_\varepsilon+\lambda\rho\le 0\right\}\cap\oli B
$$
are respectively strictly concave and strictly convex and contain $z^0$ as an interior point. Set $\oli B_3=\oli B_1\cap\oli B_2$ and $\oli D_1=\oli C^+_\l\cup \oli B_3$; $z^0$ is an interior point of $\oli D_1$ and ${\rm b}\oli B_1\smallsetminus{\rm b}\oli B_2\Subset B$. By construction $\oli{D}\cap\oli D_1\cap\oli B=\oli{D}\cap\oli B$ and $D\cap B$ is an intersection of two strictly convex domains with smooth boundaries thus applying Proposition \ref{LTL} we obtain
$$
H^r(\oli{D}\cap\oli B,\mathcal O)=\{0\}
$$
for $1\le r\le n-2$ and consequently, in view of (\ref{Hom}), the vanishing
\begin{equation}\label{BUMP2}
H^r(\oli{D}\cap\oli B,\mathcal F)=\{0\}.
\end{equation}
Iterating this procedure  we get the conclusion.
\enddemo
\bp\label{ONTO1}
For every $\varepsilon\in(0,\varepsilon_0)$ there exists $\varepsilon'<\varepsilon$ such that the homomorphism 
$$
H^r(\oli C^+_{\varepsilon'},\mathcal F)\longrightarrow H^r(\oli C^+_\varepsilon,\mathcal F)
$$
is onto for $1\le r\le p(\oli C^+_\varepsilon;\mathcal F)-2$.
\ep
\demo
Keeping the notations of Lemma \ref{BUMP1} we apply the Mayer-Vietoris exact sequence for closed sets to $\oli D_1=\oli D\cup\left(\oli D_1\cap\oli B\right)$. We get 
$$
\cdots \to H^r(\oli{D}_1,\mathcal F)\to H^r(\oli{D},\mathcal F)\oplus H^r(\oli{D}_1\cap\oli B,\mathcal F)\to H^r(\oli{D}\cap\oli D_1\cap\oli B,\mathcal F)\to \cdots\nonumber 
$$
thus in view of (\ref{BUMP2}) the homomorphism
$$ 
H^r(\oli{D}_1,\mathcal F)\to H^r(\oli{D},\mathcal F) 
$$
is onto for $1\le r\le n-2$. By induction, we obtain that  the homomorphism
$$ 
H^r(\oli{D}_m,\mathcal F)\to H^r(\oli{D},\mathcal F) 
$$
is onto for $1\le r\le p(\oli C^+_\varepsilon;\mathcal F)-2$. Since $\oli C^+_\varepsilon\subset D_m$ if $\varepsilon'<\varepsilon$ is near $\varepsilon$ one has $\oli C^+_{\varepsilon}\subset C^+_{\varepsilon'}\Subset D_m$, whence the homomorphism
$$ 
H^r(\oli C^+_{\varepsilon'},\mathcal F)\to H^r (\oli C^+_{\varepsilon},\mathcal F) 
$$
is onto for $1\le r\le p(\oli C^+_\varepsilon;\mathcal F)-2$. In particular, the canonical homomorphism 
\begin{equation}\label{ONTO2} 
H^r(C^+_{\varepsilon'},\mathcal F)\stackrel{\delta}{\to} H^r(\oli C^+_{\varepsilon},\mathcal F) 
\end{equation}
is onto for $1\le r\le p(\oli C^+_\varepsilon;\mathcal F)-2$. 
\enddemo
From Proposition \ref{ONTO1} we derive 
\bp\label{ONTO11}
For every $\varepsilon\in (0,\varepsilon_0)$ there exists an $\oli\varepsilon<\varepsilon$ such that for every $\varepsilon'\in [\oli\varepsilon,\varepsilon)$ the homomorphism
\begin{equation}\label{ONTO3}
H^r(C^+_{\varepsilon'},\mathcal F)\stackrel{\delta}{\to} H^r(\oli C^+_{\varepsilon},\mathcal F) 
\end{equation}
is onto for $1\le r\le p(\oli C^+_\varepsilon;\mathcal F)-2$.
\ep
\demo
We fix $\varepsilon_0$ as in Lemma \ref{BUMP1}. Let $\Lambda$ be the (non-empty) set of the positive numbers $\varepsilon'<\varepsilon$ such that the homomorphism (\ref{ONTO2}) is onto and  $\oli\varepsilon=\inf\,\Lambda$. It follows (cfr.\ \cite[Lemma pag. 241]{AG} for closed subsets) that the homomorphism (\ref{ONTO3}) is onto.
\enddemo
A second consequence of Proposition \ref{ONTO1} is the following finiteness theorem
\bt\label{FIN1} Under the conditions of Theorem \ref{ISO}, there exists $\varepsilon_1\le \varepsilon_0$ such that 
 $$
\dim_\C\,H^1(\oli C^+_\varepsilon,\mathcal F)<+\infty.
$$
for every $\varepsilon\in(0,\varepsilon_1)$.  
\et
\demo
We first observe the following. Let $\Omega\subset\C^n$ be a domain, $K\subset \Omega$ a compact subset. It is known that $\mathcal F(\Omega)$ is a Fr\'echet space. The space $\mathcal F(K)$ is an $\mathcal L\mathcal F$-space i.e.\ a direct limit of Fr\'echet spaces and its topology is complete (cfr.\ \cite[pag. 315]{GRO}). Moreover, the restriction 
$$
\mathcal F(\Omega)\stackrel{\d}{\longrightarrow} \mathcal F(K)
$$
is a compact map i.e.\ there exists a neighbourhood $U$ of the origin in $ \mathcal F(\Omega)$ such that $\oli {\d(U)}$ is a compact subset of $ \mathcal F(K)$. This is a consequence of the following well known fact: if $\Omega'$ is a relatively compact subdomain of $\Omega$ then the restriction $ \mathcal F(\O)\to \mathcal F(\Omega')$ is a compact map. Take $\varepsilon_0$ as in Lemma \ref{BUMP1}. The proof is similar to that of Th\'eor\`eme 11 in \cite{AG} taking into account the following facts:
\begin{itemize}
\item[1)] Leray theorem for acyclic closed coverings (see Th\'eor\`eme 5.2.4 and Corollaire in \cite{Go58})\\
\item[2)] the theorem of L. Schwartz on compact perturbations $u+v$ of a surjective linear operator $u:E\to F$ where $E$ is a Fr\'echet (see \cite[Corollaire 1]{GRO}).
\end{itemize}
\enddemo
We remark that, up to some modifications in the technical deteils of the proof, the finiteness result holds for all cohomology groups:
\bt Under the conditions of Theorem \ref{ISO} there exists $\varepsilon_1\le \varepsilon_0$ such that 
 $$
\dim_\C\,H^r(\oli C^+_\varepsilon,\mathcal F)<+\infty,
$$
for every $\varepsilon\in(0,\varepsilon_1)$ and $1\le r\le p(\oli C^+_\varepsilon;\mathcal F)-2$.  
\et
\subsection{Approximation}\label{FIN}
This subsection is devoted to approximation by global sections.
\bl\label{ONTO4}
Let ${\rm depth}(\mathcal F_z)\ge 4$ for every $z\in\{\varphi=\varepsilon\}$, $\varepsilon\in(0,\varepsilon_0)$. Then, for every $z^0\in {\rm b}C^+_\varepsilon$ there exists a closed neighbourood $\oli U$ of $z^0$ such that the homomorphism
$$
H^0(\oli U,\mathcal F)\longrightarrow H^0(\oli U\cap{ \oli C^+_\varepsilon },\mathcal F)
$$
is dense image.
\el
\demo
This is known if $z^0\in\Gamma_1\cup\Gamma_2$ i.e.\ when $z^0$ is a point of pseudoconvexity or pseudoconcavity (see \cite{AG}), thus we may assume that $z^0\in\oli\Gamma_1\cap\oli\Gamma_2$. First we consider the case $\mathcal F=\mathcal O$. We may suppose that there exists a sufficiently small closed ball $\oli B$  of positive radius, centered at $z^0$ such that $\oli B\cap \{\varphi\ge\varepsilon\}$ and $\oli B\cap \{P_\varepsilon\le 0\}$ respectively are strictly concave and strictly convex (again, locally, up to a biholomorphism). Take a real hyperplane with equation $l=0$ such that $z^0\in \{l>0\}$ and 
$\{l=0\}\cap \{\varphi\le\varepsilon\}\Subset B$. Let $\psi=\alpha\varphi-\varepsilon+\beta l$, $\alpha,\beta$ positive real numbers; $\psi$ is strongly plurisubharmonic. For $\alpha,\beta$ sufficiently small the hypersurface $\{\psi=0\}\cap\{l<0\}$ is a portion of a compact smooth hypersupersurface which bounds a domain $D\Subset B$. Set 
$$
\oli V=\{P_\varepsilon\le 0\}\cap\oli D, \>\oli W=\oli D\smallsetminus\{\varphi<\varepsilon\}
$$ 
and $\oli U'=\oli V\cap \oli W$. We are going to prove that $H^1(\oli V\cup \oli W,\mathcal O)=0$. Let $R=\oli D\smallsetminus\oli V\cup \oli W$. Since $\oli D$ is a Stein compact, from the exact sequence of cohomology relative to the closed subspace $\oli V\cup \oli W$ we get the isomorphism
\begin{equation}\label{ISO3}
H^r(\oli V\cup \oli W,\mathcal O)\simeq H^{r+1}_c(R,\mathcal O).
\end{equation}
for $r\le n-2$.
$R$ is an open subset of $S=\oli D\cap\{\varphi<\varepsilon\}$. Set $R'=S\smallsetminus R$. Again, by the cohomology with compact supports relative to the closed subspace $R'=S\smallsetminus R$ we get the exact sequence of groups   
$$
\cdots\longrightarrow H^{r}_c(S,\mathcal O)\longrightarrow  H^{r}_c(R',\mathcal O)\longrightarrow
$$
$$
\longrightarrow H^{r+1}_c(R,\mathcal O)\longrightarrow H^{r+1}_c(S,\mathcal O)\longrightarrow\cdots\>.
$$
Since $S$ and $R'$ have a fundamental system of Stein neighbourhoods (see \cite{T}) and $n\ge 3$, we have
$$
H^r_c(S,\mathcal O)=H^r_c(R',\mathcal O)=0
$$
for $1\le r\le n-2$ and consequently $H^r_c(R,\mathcal O)=0$ for $1\le r\le n-2$. In view of the isomorphism (\ref{ISO3}) we obtain
$$
H^r(\oli V\cup\oli W,\mathcal O)=0
$$
for $1\le r\le n-2$. In particular, since $n\ge 3$, (\ref{ISO3}) implies that 
$$
H^1(\oli V\cup \oli W,\mathcal O)=0,
$$
thus that every function $f\in\mathcal O(\oli U')$ is a difference of two functions $f_1-f_2$ where $f_1\in\mathcal O(\oli V)$, $f_2\in\mathcal O(\oli W)$. Since $\oli V$ is Runge in $\oli D$ there exists a sequence of holomorphic functions $f_\nu\in\mathcal O(\oli D)$ such that $f_\nu\to f_1$ in 
$\mathcal O(\oli V)$. Moreover, by the extension theorem in \cite{LT} the function $f_2$ extends holomorphically to $W\cap \{l\le 0\}$.  Choose a smooth function $\varrho\in C^{\infty}_0(D)$ such that $\varrho\ge 0$ and $\varrho (z^0)\neq 0$ and
a positive number $\lambda$ such that the closed domains
$$
\oli D_1=\left\{\varphi-{\varepsilon}+\l\varrho \le 0\right\}\cap\oli D,\>\oli D_2=\left\{P_\varepsilon-\l\rho\le 0\right\}\cap\oli D
$$
are respectively strongly pseudoconcave and strongly pseudoconvex, both contain $z^0$ as an interior point, ${\rm b}D_1\smallsetminus\{\varphi=\varepsilon\}\cap D$ is relatively compact in $D\cap \{l>0\}$ and ${\rm b}D_2\smallsetminus\{P_\varepsilon=0\}$ is relatively compact in $D$. Then we define $\oli U={\oli D}_1\cap{\oli D}_2$.

Observe that, by construction, Proposition \ref{LTL} applies, thus $H^r(\oli U\cap C^+_\varepsilon,\mathcal O)=0$ for $1\le r\le n-2$.

In the general case, since $\oli D$ is Stein, we have on $\oli D$ an exact sequence 
$$
\xymatrix{0\ar[r]&\mathcal H\ar[r]^\a&\mathcal O^q\ar[r]^\b &\mathcal F\ar[r]& 0}.
$$
Consider the following commutative diagram of continuous maps
$$
\xymatrix{H^0(\oli U, \mathcal O^q)\ar[r]^{\alpha}\ar[d]_{\sf r}& H^0(\oli U,\mathcal F)\ar[d]_{\sf r} \\  H^0(\oli U\cap \oli C^+_\varepsilon,\mathcal O^q)\ar[r]^{\beta} &H^1(\oli U\cap \oli C^+_\varepsilon,\mathcal F)}
$$
where $\sf r$ denotes the natural restriction. Then, since ${\rm depth}(\mathcal F_z)\ge 4$ for every $z\in D$, we have ${\rm depth}(\mathcal H_z)\ge 5$ for every $z\in D$. Again by Proposition 7.5 in \cite{LTL} we have $H^1(\oli U\cap C^+_\varepsilon,\mathcal F)=0$ whence the homomorphism 
$$
H^0(\oli U\cap \oli C^+_\varepsilon, \mathcal O^q)\longrightarrow H^0(\oli U\cap \oli C^+_\varepsilon, \mathcal F)
$$
is onto. Let $\sigma\in H^0(\oli U\cap \oli C^+_\varepsilon, \mathcal F)$ and $N$ a neighbourhood of $\sigma$. Let $g\in H^0(\oli U, \mathcal O^q)$ such that $\beta(g)=\sigma$. Since the homorphism
$$
H^0(\oli U, \mathcal O^q)\longrightarrow H^0(\oli U\cap \oli C^+_\varepsilon, \mathcal O^q)
$$ 
is dense image there exists $h\in H^0(\oli U, \mathcal O^q)$ such that ${\sf r}(h)\in \beta^{-1}(N)$. Then ${\sf r}(\alpha(h))\in N$ with $\alpha(h)\in H^0(\oli U,\mathcal F)$. This shows that the homomorphism
$$
H^0(\oli U, \mathcal O^q)\longrightarrow H^0(\oli U\cap \oli C^+_\varepsilon, \mathcal O^q)
$$
is dense image.
\enddemo
\bl\label{ONTO5}
Let $\mathcal F$ and $\varepsilon_0$ be as in Lemma \ref{ONTO4}. Then for every $\varepsilon\in(0,\varepsilon_0)$ there exists $\varepsilon_2\in (0,\varepsilon)$ such that for every $\varepsilon'\in (\varepsilon_2,\varepsilon)$ the homomorphism 
$$
H^0(\oli C^+_{\varepsilon'}, \mathcal F)\longrightarrow H^0(\oli C^+_\varepsilon, \mathcal F)
$$
is dense image.
\el
\demo
With the notations of Lemma \ref{BUMP1} we have 
$$
\oli D=\oli C^+_\varepsilon,\> \oli D_1=\oli D\cup \oli B,\>\oli D_1=\oli D\cup(\oli D_1\cap\oli B)
$$ 
and we set $\oli V=\oli D_1\cap\oli B$. In view of Lemma \ref{ONTO4} we may assume that the homomorphism
$$
H^0(\oli V,\mathcal F)\longrightarrow H^0(\oli V\cap \oli D,\mathcal F)
$$
is dense image. Moreover, $H^1(\oli V,\mathcal F)=0$. Let $\oli{\mathcal U}$ be the closed covering $\left\{\oli D,\oli V\right\}$ of $\oli D_1$, ${\sf Z}^1(\oli{\mathcal U},\mathcal F)$ and ${\sf B}^1(\oli{\mathcal U},\mathcal F)$ respectively the space of cocycles and coboundaries of $\oli{\mathcal U}$ with values in $\mathcal F$. Since $H^1(\oli{\mathcal U},\mathcal F)$ is a subgroup of 
$H^1(\oli D_1,\mathcal F)$ which is of finite dimension (cfr.\ Theorem \ref{FIN1}) we have
$$
\dim_\C\,H^1(\oli{\mathcal U},\mathcal F)<+\infty.
$$ 
It follows that $H^1(\oli{\mathcal U},\mathcal F)$ is of finite dimension in the $\mathcal L\mathcal F$-space ${\sf Z}^1(\oli{\mathcal U},\mathcal F)$, thus an 
$\mathcal L\mathcal F$-space for the induced topology. Moreover, in view of the Banach open mapping theorem the suriective map 
$$
H^0(\oli D,\mathcal F)\oplus H^0(\oli V,\mathcal F)\longrightarrow{\sf B}^1(\oli{\mathcal U},\mathcal F)
$$
given by $s\oplus\sigma\mapsto s_{\vert\oli D\cap\oli V}-\sigma_{\vert\oli D\cap\oli V}$ is a topological homomorphism. 

Let $s\in H^0(\oli D,\mathcal F)$; $s_{|\oli V\cap\oli D}\in {\sf B}^1(\oli{\mathcal U},\mathcal F)$. By Lemma \ref{ONTO4}, there exists a generalized sequence
 $\{s_\nu\}\subset H^0(\oli V,\mathcal F)$ such that 
 $$
 s_\nu{_{|\oli V\cap\oli D}}-s_{|\oli V\cap\oli D}\longrightarrow 0.
 $$
 In view of Banach theorem there exist two generalized sequences 
 $\sigma^1_\nu\in H^0(\oli D,\mathcal F)$, $\sigma^2_\nu\in H^0(\oli V,\mathcal F)$ such that
 $$
 {\sigma^1_\nu}_{\vert\oli D\cap\oli V}-{\sigma^2_\nu}_{\vert\oli D\cap\oli V}={s_\nu}_{\vert\oli D\cap\oli V}-{s}_{\vert\oli D\cap\oli V},
 $$
$$
 \sigma^1_\nu\to 0,\>\>\sigma^2_\nu\to 0.
 $$    
It follows that for every $\nu$
$$
\tilde s_\nu=\begin{cases}s-\sigma^1_\nu & {\rm on}\ \oli D \\s_\nu-\sigma^2_\nu &{\rm on}\ \oli V\end{cases}
$$
is a section of $\mathcal F$ on $\oli D_1$ and that $\tilde s_\nu\to s$. In order to ends the proof we apply this procedure a finite numbers of times. 
\enddemo
As a corollary we get the following
\bp\label{ONTO51}
Let $\mathcal F$ and $\varepsilon_0$ be as in Theorem \ref{ISO}. Then for every $\varepsilon\in(0,\varepsilon_0)$ there exists $\oli\varepsilon_0\in[0,\varepsilon)$ such that for every $\varepsilon'\in(\oli\varepsilon_0,\varepsilon]$ the homomorphism 
$$
H^0(C^+_{\varepsilon'}, \mathcal F)\longrightarrow H^0(\oli C^+_\varepsilon, \mathcal F)
$$
is dense image.
\ep
\demo
Let $I\subset (0,\varepsilon_0)$ be the (non-empty) set of $\varepsilon'<\varepsilon$ such that the homomorphism
$$
H^0(\oli C^+_{\varepsilon'}, \mathcal F)\longrightarrow H^0(\oli C^+_\varepsilon, \mathcal F)
$$
is dense image. Let $\oli\varepsilon=\inf\,I$ and $\left\{\varepsilon_\nu\right\}$ be a decreasing sequence with $\varepsilon_0=\varepsilon$, $\varepsilon_\nu\to \oli\varepsilon$ and set $F_\nu=H^0(\oli C^+_{\varepsilon_\nu}, \mathcal F)$. The topology of $F_\nu$ can be defined by an increasing sequence $\{p^{(\nu)}_ j\}_{j\in\N}$ of translation invariant seminorms. Let for $\nu\ge 1$
$$
{\sf r}_\nu:F_\nu\longrightarrow F_{\nu-1}
$$
be the restriction map; then
$$
H^0(\oli C^+_{\oli\varepsilon}, \mathcal F)=\varprojlim\limits_{\{\sf r_\nu\}}F_\nu
$$
and denote $\pi_\nu:H^0(\oli C^+_\varepsilon, \mathcal F)\to F_\nu$ the natural map. We have to show that $\pi_0$ is dense image.

 Let $s\in F_0=H^0(\oli C^+_\varepsilon, \mathcal F)$ and $N$ a neighbourhood of $s_0$. We may assume that 
$$
N=\left\{s\in F_0:p^{(0)}_0(s-s_0)<\varepsilon\right\}.
$$
Since the maps $\sf r_\nu$ are continuous and dense image we can choose elements $s_\nu\in F_\nu$, for $\nu\ge 0$, satisfying the following conditions:
\begin{itemize}
\item[] $s_1\in F_1\>\>\>\> p^{(0)}_0\left({\sf r}_1(s_1)-s_0\right)<\varepsilon/2 $
\item[]$ s_2\in F_2\>\>\>\>  p^{(1)}_0\left({\sf r}_2(s_2)-s_1\right)<\varepsilon/2$
\begin{itemize}
\item[]$\>\>\>\>\>\>\>\>\> p^{(0)}_1\left({\sf r}_1{\sf r}_2(s_2)-{\sf r}_1(s_1)\right)<\varepsilon/2^2$
\end{itemize}
\item[]$s_3\in F_3\>\>\>\>  p^{(2)}_0\left({\sf r}_3(s_3)-s_2\right)<\varepsilon/2$
\begin{itemize}
\item[]$\>\>\>\>\>\>\>\>\> p^{(1)}_1\left({\sf r}_2{\sf r}_3(s_3)-{\sf r}_2(s_2)\right)<\varepsilon/2^2$
\item[]$\>\>\>\>\>\>\>\>\> p^{(0)}_2\left({\sf r}_1{\sf r}_2{\sf r}_3(s_3)-{\sf r}_1{\sf r}_2(s_2)\right)<\varepsilon/2^3$
\end{itemize}
\end{itemize}
and so on. Then, for every $\nu\in\N$, the series
$$
s_\nu+\left({\sf r}_{\nu+1}(s_{\nu+1})-s_\nu\right)+\left({\sf r}_{\nu+1}{\sf r}_{\nu+2}(s_{\nu+2})-{\sf r}_{\nu+1}(s_{\nu+1})\right)+\ldots
$$
is convergent in $F_\nu$ and ${\sf r}_{\nu}(\sigma_{\nu})=\sigma_{\nu-1}$. Hence $\sigma=\{\sigma_\nu\}_{\nu\in\N}$ belongs to $H^0(C^+_{\oli\varepsilon},\mathcal F)$ and, by definition  $p^{(0)}_0(\sigma_0-s_0)<\varepsilon$, i.e.\ $\pi_0(\sigma_0)\in N$.
\enddemo

\noindent{\bf Proof of Theorem \ref{ISO}.} The proof uses Corollary \ref{ONTO1} and Lemma \ref{ONTO4}. With the notations of Lemma \ref{BUMP1} we have 
$$
\oli D=\oli C^+_\varepsilon,\> \oli D_1=\oli D\cup \oli B,\>\oli D_1=\oli D\cup(\oli D_1\cap\oli B),\>\oli D\cap(\oli D_1\cap\oli B)=\oli D\cap\oli B.
$$ 
We may assume that the homomorphism
$$
H^1(\oli D_1,\mathcal F)\longrightarrow H^1(\oli D,\mathcal F)
$$
is onto and
$$
H^0(\oli D_1,\mathcal F)\longrightarrow H^0(\oli D,\mathcal F)
$$
is dense image. Moreover, $H^1(\oli B\cap\oli D_1,\mathcal F)=0$. Thus it is sufficient to show that the homomorphism
$$
H^1(\oli D_1,\mathcal F)\longrightarrow H^1(\oli D,\mathcal F)
$$
is injective. 

Since  $H^1(\oli D_1\cap\oli B,\mathcal F)=0$ the Mayer-Vietoris exact sequence applied to $\oli D_1=\oli D\cup(\oli D_1\cap\oli B)$ gives the exact sequence
$$
\xymatrix{H^0(\oli D\cap\oli B,\mathcal F)\ar[r]^{\sf a}& H^1(\oli D_1,\mathcal F)\ar[r]^{\sf b}& H^1(\oli D,\mathcal F)&}.
$$
Let $\xi\in {\sf Ker}\,\sf b={\sf Im}\,\sf a$, $\xi={\sf a}(\eta)$ with $\eta\in H^0(\oli D\cap\oli B,\mathcal F)$.  By Lemma \ref{ONTO4} $\eta$ is approximated by a sequence  $\left\{\eta_\nu\right\}\subset H^0(\oli D_1\cap\oli B,\mathcal F)$. Each $\eta_\nu$ is a $1$-coboundary of the closed covering $\mathcal U=\left\{\oli D,\oli D_1\cap\oli B\right\}$ with values in $\mathcal F$ and such a space is closed in the space ${\sf Z}^1\left(\mathcal U,\mathcal F\right)$ of the $1$-cocycles. This proves that $\eta$ is a $1$-coboundary of $\left\{\mathcal U,\mathcal F\right\}$, whence  $\xi={\sf a}(\eta)=0$.
\enddemo

\br\label{JUMP} 
In the full $q$-corona the cohomology of all coronae are isomorphic (see \cite{AG}). Differently, in the semi $1$-corona case the cohomology groups are isomorphic up to a critical $\oli\varepsilon$, where the dimension of the cohomology spaces jumps, then they are again all isomorphic up to a second critical value, and so on. They must not be all isomorphic, nor they dimensions must be bounded.
\er

\section{Extension of coherent sheaves and analytic subsets}
An interesting consequence is that on a semi  $1$-corona $C^+=C^+_{0,1}$ Theorem A of Oka-Cartan-Serre holds for a coherent sheaf $\mathcal F$ satisfying the conditions of Theorem \ref{ISO}. We first prove the following
\bl\label{ISO5}
Let $X$ be a complex space, $\mathcal F\in{\sf Coh}(X)$ satisfying the following property: for every $x\in X$ there exists a subset $Y\not\ni x$ of $X$ such that:
\begin{itemize}
\item[i)]
$
H^1(X,\mathcal F)\simeq H^1(Y,\mathcal F)
$
\vspace{2mm}
\item[ii)] if $\mathcal M_{[x]}$ denotes the ideal of $\{x\}$ the homomorphism
$$
H^1(X,\mathcal M_{[x]}\mathcal F)\longrightarrow H^1(Y,\mathcal M_{[x]},\mathcal F) 
$$
is injective.
\end{itemize} 
Then, for every $x\in X$ the space $H^0(X,\mathcal F)$ of the global sections of $\mathcal F$ generates $\mathcal F_x$ over $\mathcal O_{X,x}$.
\el
\demo
Let $x\in X$ and $Y$ satisfying the conditions of the lemma. Consider the exact sequence of sheaves
$$
\xymatrix{0\ar[r]& \mathcal M_{[x]}\mathcal F \ar[r]& \mathcal F\ar[r]{}& \mathcal F/\mathcal M_{[x]} \mathcal F\ar[r]^{}&0}
$$
and the associated diagram
$$
H^0(X,\mathcal F)\to H^0(X, \mathcal F/\mathcal M_{[x]} \mathcal F)\to\xymatrix{H^1(X, \mathcal M_{[x]}\mathcal F)\ar[r]^{\delta}\ar[d]_{\alpha} & H^1(X,\mathcal F)\ar[d]_{\beta} \\  H^1(Y,\mathcal M_{[x]}\mathcal F)\ar[r]^{\gamma} &H^1(Y,\mathcal F).}
$$
The homomorphism is injective by hypothesis and $\beta$ is an isomorphism since $\mathcal M_{{[x]}\vert Y}\simeq \mathcal F_{\vert Y}$, thus $\gamma$ is an isomorphism. It follows that $\delta$ is injective and consequently that the homomorphism
$$
H^0(X, \mathcal F)\to H^0(X, \mathcal F/ \mathcal M_{[x]}\mathcal F)\simeq \mathcal F_x/\mathcal M_{{[x]},x}\mathcal F_x
$$
is onto. Then the Lemma of Nakayama implies that 
$$
H^0(X, \mathcal F)\longrightarrow\mathcal F_x
$$
is onto and this proves the lemma.
\enddemo
Keeping the notations of the proof of Theorem \ref{ISO}, we deduce the following
\bc\label{OKS1}
Under the conditions of Theorem \ref{ISO} for every compact subset 
$$
K\subset C^+_{\varepsilon'}\cap\left\{\varphi>\varepsilon'\right\}\cap\left\{P_{\varepsilon'}<0\right\}
$$
there exist sections $s_1,\ldots,s_k\in H^0( C^+_{\varepsilon'},\mathcal F)$ which generate $\mathcal F_z$ for every $z\in K$.
\ec

\bt\label{OKS}
Let $C^+=\left(B_{1}\smallsetminus B_0\right)\cap\{h\ge 0\}$ and $\mathcal F\in {\sf Coh}(C^+)$. If ${\rm depth}(\mathcal F_z)\ge 3$ on $\{\varphi=0\}$ then for every $a>0$ near $0$ $\mathcal F_{|B_1\smallsetminus\oli B_a}$ extends on $B_1\cap\{h\ge 0\}$ by a coherent sheaf $\widetilde{\mathcal F}_a$.
\et
\demo
With the usual notations choose $\varepsilon_0\in (0,a)$, and $c_0>0$ such that 
\begin{itemize}
\item[i)] $\mathcal F$ is defined on the semi $1$-corona $\left(B_1\smallsetminus B_{-\varepsilon}\right)\cap \left\{h>-c\right\}$
\item[ii)] $\left\{z\in B_1:h(z)\ge c\right\}\Subset \left\{z\in B_1:P_\varepsilon(z)<0\right\}$
\item[iii)] for every $\varepsilon\in(0,\varepsilon_0)$, $c\in(0,c_0)$ the hypersurfaces $\left\{P_\varepsilon=-c\right\}$, $P_\varepsilon(z)=\varepsilon\vert z\vert^2-h$, meet the hypersurfaces 
$\{\varphi=\varepsilon\}$, $\{\varphi=-\varepsilon\}$ transversally. 
\end{itemize}
Let $Y^+_{\a,\b}$ denote the semi $1$-corona $\left\{\alpha<\varphi<\beta\right\}\cap\left\{h>c\right\}$, with $\a<\b<\varepsilon$. In view of Corollary \ref{OKS1} applied to the semi $1$-corona $Y^+_{\varepsilon,a}$ there exist $\a,\b,\gamma\in(0,a)$ with $\a<\b<\gamma$ such that $H^0(Y^+_{\a,\gamma},\mathcal F)$ generates $\mathcal F$ on $K_{\b,\gamma}=\oli Y^+_{\b,\gamma}\cap\left\{h\ge 0\right\}$. Thus on $K_{\b,\gamma}$ there exists an exact sequence
$$
\xymatrix{\mathcal O^p\ar[r]^\beta &\mathcal F\ar[r]& 0}.
$$
Since, by hypothesis, ${\rm depth}(\mathcal F_z)\ge 3$ for every $z\in K_{\b,\gamma}$ we have 
$$
{\rm depth}({\sf Ker}\,\alpha)\ge 4
$$
on $K_{\b,\gamma}$ (cfr.\ \cite{CART1}). Again by Corollary \ref{OKS1} there exist $\b_1, \gamma_1\in (\b,\gamma)$, $\b_1<\gamma_1$ and sections $\sigma_1,\ldots,\sigma_l$ on $ K_{\b_1,\gamma_1}=\oli Y^+_{\b_1,\gamma_1}\cap\left\{h\ge 0\right\}$ which generate $({\sf Ker}\,\alpha)_z$ for every $z\in V$. Since ${\sf Ker}\,\alpha$ is a subsheaf of $\mathcal O^p$, by the theorem in \cite{LT} the sections $\s_1,\ldots,\s_l$ extend holomorphically on 
$$
\left\{\varphi\le\gamma_1\right\}\cap\left \{h\ge 0\right\}
$$
and their extensions $\tilde \s_1,\ldots,\tilde \s_l$ generate a coherent sheaf $\mathcal H$ on 
 $$
 \left\{\varphi\le\gamma_1\right\}\cap\left \{h\ge 0\right\}.
$$
Let $\widetilde{\mathcal F'}_a$ be the sheaf defined by
$$
\widetilde{\mathcal F'}_{a,z}=\begin{cases}\mathcal F_z& {\rm for}\ z\in\left\{\varphi>\gamma_1\right\}\cap\left \{h\ge 0\right\}
\\ \mathcal O_z/{\mathcal H}_z &{\rm for}\ z\in\left\{\varphi\le\gamma_1\right\}\cap\left \{h\ge 0\right\};\end{cases}
$$
 $\widetilde{\mathcal F'}_{\varepsilon}$ is a coherent sheaf on $B^+_{c}\cap\{y_n>\varepsilon\}$ extending $\mathcal F$.
\enddemo
\bc\label{OKA2}
Let $X^+=\left(B_1\smallsetminus \oli B_0\right)\cap\left\{h> 0\right\}$ and $Y$ an analytic subset of $X^+$ such that ${\rm depth}(\mathcal O_{Y,z})\ge 3$ for $z$ near $\left\{\varphi=0\right\}$. Then $Y$ extends on $B_1\cap\{h\ge 0\}$ by an analytic subset. 
\ec
\demo
We apply Theorem \ref{OKS} to $X^+\cap\left\{h\ge\varepsilon\right\}$, where $\varepsilon\sim 0$ is positive. Then, for $\nu\in\N$ there  exists a coherent sheaf $\widetilde{\mathcal O}^{(\nu)}_Y$ on $B_1\cap\left\{h\ge 0\right\}$ which extends $\mathcal O_Y$; $\widetilde {Y}^{(\nu)}={\rm supp}\>\>
\widetilde{\mathcal O}^{(\nu)}_Y$ is an analytic subset extending $Y\cap \left(B_1\smallsetminus B_{1/\nu}\right)\cap\left\{h\ge\varepsilon\right\}$. In view of the strong pseudoconvexity of ${\rm b}B_{1/\nu}$, the subset $F_\nu=\widetilde {Y}^{(\nu)}\smallsetminus\widetilde {Y}^{(\nu+1)}$ is a finite set of points which is contained in $B_{1/\nu}$. Start by $\nu=2$ and consider the first extension $\widetilde {Y}^{(2)}$. Then $\widetilde {Y}^{(2)}\smallsetminus F_2\cap\left (B_{1/2}\smallsetminus B_{1/3}\right)$ coincide with $Y$ on $\left (B_1\smallsetminus B_{1/3}\right)$ and so on. To handle different extensions depending on $\varepsilon$ we argue in the same way.
\enddemo

\end{document}